\documentstyle[fleqn,12pt]{article}
\begin{document}
\pagenumbering{arabic}
\setcounter{page}{1}
\pagestyle{plain}
\baselineskip=16pt

\thispagestyle{empty}
\rightline{MSUMB 99-1, January 1999} 
\vspace{1.4cm}

\begin{center}
{\Large\bf Differential geometry of GL$_h(1\vert 1)$} 
\end{center}

\vspace{1cm}
\begin{center} Salih \c Celik 
\footnote{{\bf E-mail address}: sacelik@yildiz.edu.tr}\\
Yildiz Technical University, Department of Mathematics, \\
34210 Davutpas-Esenler, Istanbul, TURKEY.
\end{center}

\vspace{1.5cm}
{\bf Abstract}

We construct a right-invariant differential calculus on the 
quantum supergroup GL$_h(1\vert 1)$ and obtain the $h$-deformed 
superalgebra of GL$_h(1\vert 1)$.

\vfill\eject
\noindent
{\bf I. INTRODUCTION}

In the last few years, the theory of quantum (super) groups like GL(2), 
GL$(1\vert 1)$, etc., were generalized in two ways. Both of the 
generalizations are based on the deformation of the algebra of functions 
on the matrix (super) groups generated by coordinate functions $T_j^i$ 
which normally commute. These deformations of Lie (super) groups are 
algebraic structures depending on one (or more) continuous parameter. 
We have a standard Lie (super) group for particular values of the 
deformation parameters. Quantum (super) groups$^{1-3}$ present the 
examples of (graded) Hopf algebras. They have found application in 
diverse areas of physics and mathematics$^4$. 

The $q$-deformation of Lie (super) groups can be realized on a quantum 
(super) space in which coordinates are noncommuting$^2$. Recently the 
differential calculus on noncommutative (super) space has been intensively 
studied both by mathematicians and mathematical physicists. There is much 
activity in differential geometry on quantum groups. Throughout the recent 
development of differential calculus on the quantum groups two principal 
concepts are readily seen. First of them, formulated by Woronowicz$^5$, 
is known as bicovariant differential calculus on the quantum groups. 
Another concept, introduced by Woronowicz$^6$ and Schirmacher {\it et al}$^7$ 
proceeds from the requirement of a calculus only. There are many papers in 
this field$^8$. We shall consider the second concept. 

Another type of deformation, the so called $h$-deformation, which is a 
new class of quantum deformations of Lie groups and Lie algebras has recently 
been intensively studied$^9$. This deformation may be obtained as a 
contraction of the $q$-deformation$^{10}$. There is much interest in studies 
relating to various aspects of the $h$-deformed algebra. The differential 
geometry of SL$_h(2)$ was given in$^{11}$. In this work, we introduce a 
right-invariant differential calculus on the quantum supergroup 
GL$_h(1\vert 1)$. This quantum supergroup was obtained in ref. 12 using a 
contraction procedure given in Ref. 10. 

Let us briefly discuss the content of the paper. In the second 
section, the basic notations of the Hopf algebra structure on the 
quantum supergroup GL$_h(1\vert 1)$ are introduced. In the 
third section we shall obtain the commutation relations for the 
group parameters (the matrix elements) and their differentials 
so we have a differential algebra. This differential algebra (extended 
algebra) has a Hopf algebra structure. Later, we shall construct the 
Cartan-Maurer one-forms and obtain the needed commutation relations. 
Using these commutation relations, we shall describe the quantum 
superalgebra for the vector fields (superalgebra generators) 
for GL$_h(1\vert 1)$ and derive the commutation relations between 
the group parameters and the algebra generators. We shall also show that the 
obtained quantum superalgebra can be rederived using the partial derivatives 
and their relations. 

\noindent
{\bf II. THE ALGEBRA OF FUNCTIONS ON GL$_h(1\vert 1)$}

\noindent
Elementary properties of quantum supergroup GL$_h(1\vert 1)$ are 
described in Ref. 12. We state briefly the properties we are going to 
need in this work. Here we denote $q$-deformed objects by primed 
quantities. Unprimed quantities will be represent transformed 
coordinates. As usual, it known that even (bosonic) objects commute 
with everything and odd (Grassmann) objects anticommute among 
themselves. In this work, 
to obtain the quantum supergroup GL$_h(1\vert 1)$ $^{12}$, 
we shall only assume that odd elements 
$\beta$ and $\gamma$ anticommute with the 'new' deformation parameter $h$. 

Let us begin with the $q$-deformed counterparts of GL$(1\vert 1)$. The 
quantum supergroup GL$_q(1\vert 1)$ is defined by the matrices 
of the form 
$$T' = \left(\matrix{ a' & \beta' \cr \gamma' & d' \cr} \right),$$ 
where the matrix entries satisfy the following commutation relations$^{2,3}$ 
$$ a' \beta' = q \beta' a', \qquad d' \beta' = q \beta' d', $$
$$ a' \gamma' = q \gamma' a', \qquad d' \gamma' = q \gamma' d', \eqno(1)$$
$$ \beta' \gamma' + \gamma' \beta' = 0, \qquad \beta'^2 = 0 = \gamma'^2, $$
$$ a' d' = d' a' + (q - q^{-1}) \gamma' \beta'. $$

We now consider the following similarity transformation$^{10}$: 
$$T = \left(\matrix{ a & \beta \cr \gamma & d \cr} \right) = g^{-1} T' g, 
  \eqno(2) $$
where 
$$g = \left(\matrix{ 1 & 0 \cr h/(q - 1) & 1 \cr} \right), 
  \qquad h^2 = 0. \eqno(3) $$ 
Assuming $\beta$ and $\gamma$ that anticommute with the Grassmann number 
$h$ and substituting (2) into (1), we arrive at the following relations$^{12}$ 
$$a \beta = \beta a, \qquad 
  a \gamma = \gamma a + h a^2 (1 - {\cal D}_h^{-1}), $$
$$d \beta = \beta d, \qquad d \gamma = \gamma d + h d^2 ({\cal D}_h - 1), $$
$$ \beta^2 = 0, \qquad \gamma^2 = h \gamma d (1 - {\cal D}_h), \eqno(4)$$
$$ \beta \gamma = - \gamma \beta + h \beta d (1 - {\cal D}_h), $$
$$ a d = d a + h \beta d ({\cal D}_h - 1), $$
where 
$${\cal D}_h = a d^{-1} - \beta d^{-1} \gamma d^{-1} \eqno(5)$$
is the quantum superdeterminant of $T$. It can be checked that ${\cal D}_h$ 
commutes with all matrix elements of $T$. Note that, by imposing the relation 
$${\cal D}_h = 1 $$
as an interesting case, we obtain the classical special supergroup 
SL$(1\vert 1)$, instead of SL$_h(1\vert 1)$. In other words, the restriction 
of the superdeterminant  to unity does not give the quantum supergroup 
SL$_h(1\vert 1)$. It known that, in the $q$-deformed case, 
the restriction to unity 
(${\cal D}_q = a d^{-1} - \beta d^{-1} \gamma d^{-1} = 1$) 
gives the quantum supergroup SL$_q(1\vert 1)$. 

Let us denote the algebra generated by the elements $a$, $\beta$, $\gamma$, 
$d$ with the relations (4) by ${\cal A}$. We know that the algebra 
${\cal A}$ is a graded Hopf algebra with the 
following co-structures: the usual coproduct 
$$\Delta : {\cal A} \longrightarrow {\cal A} \otimes {\cal A}, \qquad 
   \Delta(T^i_j) = T^i_k \otimes T^k_j,  \eqno(6)$$
the counit 
$$ \varepsilon : {\cal A} \longrightarrow {\cal C}, \qquad 
    \varepsilon(T^i_j) = \delta^i_j, \eqno(7)$$
and the coinverse $S : {\cal A} \longrightarrow {\cal A}$, 
$$S(T) = T^{-1} = 
  \left(\matrix{ 
    a^{-1} + a^{-1} \beta d^{-1} \gamma a^{-1} & - a^{-1} \beta d^{-1} \cr 
    - d^{-1} \gamma a^{-1} & d^{-1} + d^{-1} \gamma a^{-1} \beta d^{-1} \cr }
     \right).  \eqno(8)$$
It is not difficult to verify the following properties of the 
co-structures: 
$$(\Delta \otimes \mbox{id}) \circ \Delta = 
  (\mbox{id} \otimes \Delta) \circ \Delta, \eqno(9\mbox{a})$$
$$\mu \circ (\varepsilon \otimes \mbox{id}) \circ \Delta 
  = \mu' \circ (\mbox{id} \otimes \varepsilon) \circ \Delta, \eqno(9\mbox{b})$$
$$m \circ (S \otimes \mbox{id}) \circ \Delta = \varepsilon 
  = m \circ (\mbox{id} \otimes S) \circ \Delta, \eqno(9\mbox{c})$$
where id denotes the identity mapping, 
$\mu : {\cal C} \otimes {\cal A} \longrightarrow {\cal A},~ 
  \mu' : {\cal A} \otimes {\cal C} \longrightarrow {\cal A} $ 
are the canonical isomorphisms, defined by 
$\mu(k \otimes a) = ka = \mu'(a \otimes k), ~\forall a \in {\cal A}, 
~\forall k \in {\cal C}$ 
and $m$ is the multiplication map 
$m : {\cal A} \otimes {\cal A} \longrightarrow {\cal A}, ~ 
  m(a \otimes b) = ab$. 

The multiplication in ${\cal A} \otimes {\cal A}$ follows the rule 
$$(A \otimes B) (C \otimes D) = (-1)^{p(B) p(C)} AC \otimes BD, \eqno(10)$$
where $p(X)$ is the $z_2$-grade of $X$. 

\noindent
{\bf III. DIFFERENTIAL CALCULUS ON GL$_h(1\vert 1)$}

In this section, we shall build up the right-invariant differential 
calsulus on the quantum supergroup GL$_h(1\vert 1)$. The differential 
calculus on the quantum supergroups involves functions on the supergroup, 
differentials and differential forms. 

\noindent
{\bf A. Differential algebra}

We first note that the properties of the exterior differential. 
We can introduce the exterior differential {\sf d} to be an operator 
that is nilpotent and obeys the graded Leibniz rule: 
$$ {\sf d}^2 = 0, \eqno(11\mbox{a})$$
and 
$${\sf d} (fg) = ({\sf d} f) g + (-1)^{p(f)} f ({\sf d} g), \eqno(11\mbox{b})$$
where $f$ and $g$ are functions of the group parameters. Note that, since 
the deformation parameter $h$ is an odd (Grassmann) number it must be 
anticommute with the exterior differential {\sf d}. In fact 
$$ {\sf d} (h f) = (-1)^{p(h)} h ({\sf d} f) = - h ({\sf d} f) 
   ~\Longrightarrow~ {\sf d} h = - h {\sf d}. \eqno(11\mbox{c})$$

We have seen, in the previous section, that ${\cal A}$ is an associative 
algebra generated by the matrix elements of $T$ with the relations (4). 
A differential algebra on ${\cal A}$ is a $z_2$-graded associative algebra 
$\Gamma$ equipped with an operator {\sf d} given in (11). Also 
the algebra $\Gamma$ has to be generated by ${\cal A} \cup {\sf d} {\cal A}$. 

Firstly, we wish to obtain the relations between the matrix elements of 
$T$ and their differentials. To do this, we shall use the method of 
ref. 13. For this reason, we decompose the algebra ${\cal A}$ into 
subalgebras. We denote by ${\cal A}_{a'\beta'}$ the algebra generated 
by the elements $a'$ and $\beta'$ with the relations 
$$a' \beta' = q \beta' a', \qquad \beta'^2 = 0. \eqno(12)$$
Then, a possible set of commutation relations between generators 
of ${\cal A}_{a'\beta'}$ and ${\sf d} {\cal A}_{{\sf d} a' {\sf d} \beta'}$ is 
of the form 
$$a' ~{\sf d} a' = A_1 {\sf d} a' ~a', $$
$$a' ~{\sf d} \beta' = F_{11} {\sf d} \beta' ~a' + F_{12} {\sf d} a' ~\beta', 
  \eqno(13)$$
$$\beta' ~{\sf d} a' = F_{21} {\sf d} a' ~\beta' + F_{22} {\sf d} \beta' ~a', $$
$$\beta' ~{\sf d} \beta' = A_2 {\sf d} \beta' ~\beta', $$
where the coefficients $A_i$ and $F_{ij}$ are related to the complex 
deformation parameter $q$. To determine them we use the consistency 
of calculus (see, for details, ref. 13). Continuing in this way, we 
can obtain the other relations. 

Let us now substitute the matrix elements of ${\sf d} T'$, 
$${\sf d} T' = \left(\matrix{ \alpha' & b' \cr c' & \delta' \cr} \right) = 
               \left(\matrix{\alpha - {h\over {q - 1}} b & b \cr 
       c + {h\over {q - 1}}(\delta - \alpha) & \delta - {h\over {q - 1}} b\cr} 
  \right) \eqno(14) $$
and $T'$ into (13). After rather complicated and tedious calculations by 
using the consistency of calculus, as the final result one has the 
following commutation relations 
$$a \alpha = \alpha a + h (\alpha \beta - b a), \qquad 
  a b = b a - h b \beta, $$
$$a c = c a + h (\alpha a - c \beta + \delta a), $$
$$a \delta = \delta a + h (b a + \delta \beta), $$
$$\beta \alpha = - \alpha \beta + h b \beta, \qquad 
  \beta b = b \beta, $$
$$\beta c = c \beta + h (\alpha + \delta) \beta, \qquad 
  \beta \delta = - \delta \beta - h b \beta, $$
$$\gamma \alpha = - \alpha \gamma + h (\alpha a + \alpha d + b \gamma), 
  \eqno(15)$$
$$\gamma b = b \gamma + h b (a + d), $$
$$\gamma c = c \gamma + h (\alpha \gamma + c a + c d + \delta \gamma), $$
$$\gamma \delta = - \delta \gamma + h (\delta a + \delta d - b \gamma), $$
$$d \alpha = \alpha d - h (\alpha \beta + b d), \qquad 
  d b = b d + h b \beta, $$
$$d c = c d + h (\alpha d + c \beta + \delta d), $$
$$d \delta =  \delta d + h (b d - \delta \beta). $$
It is easy to verify that the deformation parameter $h$ anticommutes with 
$\alpha$ and $\delta$. In fact, since $a h = h a$ we have 
$$0 = {\sf d} (a h - h a) = \alpha h + h \alpha.$$

To find the commutation relations between differentials, we apply the 
exterior differential {\sf d} on the relations (15) and use the nilpotency 
of {\sf d} with (11c). Then it is easy to see that 
$$\alpha b = b \alpha + h b^2, \qquad 
  \alpha c = c \alpha + h (c b + \delta \alpha), $$
$$\delta b = b \delta  - h b^2, \qquad 
  \delta c = c \delta  - h (c b - \alpha \delta), $$
$$\alpha^2 = h \alpha b, \qquad 
  \alpha \delta = - \delta \alpha + h(\delta - \alpha)b, \eqno(16)$$
$$\delta ^2 = - h \delta  b, \qquad 
  b c = c b + h(\delta + \alpha)b. $$
These relations are the relations of Gr$_h(1\vert 1)$ in ref. 14. Note that, 
the central element of ${\sf d} {\cal A}$, which is generated by the 
elements $\alpha$, $b$, $c$, $\delta$ with the relations (16), is 
$$\hat{\cal D} = b c^{-1} - \alpha c^{-1} \delta c^{-1}. \eqno(17)$$
However, the element $\hat{\cal D}$ also commutes with the generators of 
${\cal A}$. So the element $\hat{\cal D}$ is the central element of the 
algebra ${\cal A}$, too. Thus $\hat{\cal D}$ is the central element of 
the differential algebra $\Gamma$. 

An interesting note is also the following. The central element of the 
$q$-deformed differential algebra$^{13}$ is only the element $\hat{\cal D}$. 
However the superdeterminant of $T \in $GL$_h(1\vert 1)$ commutes with the 
generators of ${\sf d} {\cal A}$, too. So the superdeterminant ${\cal D}_h$ is 
also a central element for the $h$-deformed differential algebra $\Gamma$. 

\noindent
{\bf B. Hopf algebra structure on $\Gamma$}

We first note that consistency of a differential calculus with commutation 
relations (4) means that the algebra $\Gamma$ is a graded associative algebra 
generated by the elements of the set 
$\{a, \beta, \gamma, d, \alpha, b, c, \delta\}$. 
So, it is sufficent to only describe the actions of co-maps on the 
subset $\{\alpha, b, c, \delta\}$. 

We consider a map 
$\phi_R : \Gamma \longrightarrow \Gamma \otimes {\cal A}$ such that 
$$\phi_R \circ {\sf d} = ({\sf d} \otimes \mbox{id}) \circ \Delta. \eqno(18)$$
and define a map $\Delta_R$ as follows: 
$$\Delta_R(u_1 {\sf d}v_1 + {\sf d}v_2 u_2) = 
  \Delta(u_1) \phi_R({\sf d}v_1) + \phi_R({\sf d}v_2) \Delta(u_2). \eqno(19)$$
Then it can be checked that the map $\Delta_R$ leaves invariant the relations 
(15) and (16). One can also check that the following identities are satisfied: 
$$(\Delta_R \otimes \mbox{id}) \circ \Delta_R = 
  (\mbox{id} \otimes \Delta) \circ \Delta_R 
  \qquad (\mbox{id} \otimes \epsilon) \circ \Delta_R = \mbox{id}. \eqno(20)$$
However, we do not have a coproduct for the differential algebra because the 
map $\phi_R$ does not gives an analog for the derivation property (11), yet. 
So we consider another map 
$\phi_L : \Gamma \longrightarrow {\cal A} \otimes \Gamma$ 
such that 
$$\phi_L \circ {\sf d} = (\tau \otimes {\sf d}) \circ \Delta \eqno(21)$$
and a map $\Delta_L$ with again (19) by replacing $L$ with $R$. 
Here $\tau: \Gamma \longrightarrow \Gamma$ is the linear map of degree zero 
which gives $\tau(u) = (-1)^{p(u)} u$. The map $\Delta_L$ also leaves 
invariant the relations (15) and (16), and the following identities are 
satisfied: 
$$(\mbox{id} \otimes \Delta_L) \circ \Delta_L = 
  (\Delta \otimes \mbox{id}) \circ \Delta_L 
  \qquad (\epsilon \otimes \mbox{id}) \circ \Delta_L = \mbox{id}. \eqno(22)$$

To denote the coproduct, counit and coinverse which will be defined on the 
algebra $\Gamma$ with those of $\cal A$ may be inadvisable. For this reason, 
we shall denote them with a different notation. Let us define the map 
$\hat{\Delta}$ as 
$$\hat{\Delta} = \Delta_R + \Delta_L \eqno(23)$$
which will allow us to define the coproduct of the differential algebra. 
We denote the restriction of $\hat{\Delta}$ to the algebra ${\cal A}$ by 
$\Delta$ and the extension of $\Delta$ to the differential algebra $\Gamma$ by 
$\hat{\Delta}$. It is possible to interpret the relation 
$$\hat{\Delta}\vert_{\cal A} = \Delta \eqno(24)$$
as the definition of $\hat{\Delta}$ on the generators of $\cal A$ and (23) as 
the definition of $\hat{\Delta}$ on differentials. One can see that 
$\hat{\Delta}$ is a coproduct for the differential algebra $\Gamma$ where 
$$\hat{\Delta}({\sf d} T^i_j) = {\sf d} T^i_k \otimes T^k_j + 
  (-1)^{p(T^i_k)} T^i_k \otimes {\sf d} T^k_j. \eqno(25)$$
It isnot difficult to verify the following conditions: 

{\bf a)} $\Gamma$ is an ${\cal A}$-bimodule, 

{\bf b)} $\Gamma$ is an ${\cal A}$-bicomodule with left and right coactions 
$\Delta_L$ and $\Delta_R$, respectively, making $\Gamma$ a left and right 
${\cal A}$-comodule with (20) and (22), and 
$$(\Delta_L \otimes \mbox{id}) \circ \Delta_R = 
  (\mbox{id} \otimes \Delta_R) \circ \Delta_L 
\eqno(26)$$
which is the ${\cal A}$-bimodule property. So, the triple 
$(\Gamma, \Delta_L, \Delta_R)$ is a bicovariant bimodule over Hopf algebra 
${\cal A}$. In additional, since 

{\bf c)} $(\Gamma, {\sf d})$ is a first order differential calculus over 
${\cal A}$, and 

{\bf d)} {\sf d} is both a left and a right comodule map, i.e. for all 
$u \in {\cal A}$ 
$$(\tau \otimes {\sf d}) \Delta(u) = \Delta_L({\sf d} u), \qquad 
  ({\sf d} \otimes \mbox{id})\Delta(u) = \Delta_R({\sf d} u), \eqno(27)$$
the quadruple 
$(\Gamma, d, \Delta_L, \Delta_R)$ is a first order bicovariant differential 
calculus over Hopf algebra ${\cal A}$. 
 
Now let us return Hopf algebra structure of $\Gamma$. If we define a counit 
$\hat{\epsilon}$ for the differential algebra as 
$$\hat{\epsilon} \circ {\sf d} = {\sf d} \circ \epsilon = 0 \eqno(28)$$
and 
$$\hat{\epsilon}\vert_{\cal A} = \epsilon, \qquad 
  \epsilon\vert_\Gamma = \hat{\epsilon}. \eqno(29)$$
we have 
$$\hat{\epsilon}({\sf d} T^i_j) = 0, \eqno(30)$$
where 
$$\hat{\epsilon}(u_1 {\sf d}v_1 + {\sf d}v_2 u_2) = 
  \epsilon(u_1) \hat{\epsilon}({\sf d}v_1) + 
  \hat{\epsilon}({\sf d}v_2) \epsilon(u_2). \eqno(31)$$
Here we used the fact that ${\sf d}(1) = 0$. 

As the next step we obtain a coinverse $\hat{S}$. For this, it suffices to 
define $\hat{S}$ such that 
$$\hat{S} \circ {\sf d} = {\sf d} \circ S \eqno(32)$$
and 
$$\hat{S}\vert_{\cal A} = S, \qquad 
  S\vert_\Gamma = \hat{S} \eqno(33)$$
where 
$$\hat{S}(u_1 {\sf d}v_1 + {\sf d}v_2 u_2) = 
  \hat{S}({\sf d}v_1) S(u_1) + S(u_2) \hat{S}({\sf d}v_2). \eqno(34)$$
Thus the action of $\hat{S}$ on the generators $\alpha$, $b$, $c$ and 
$\delta$ is as follows: 
$$\hat{S}({\sf d} T^i_j) 
  = - (-1)^{p[(T^{-1})^i_k]} (T^{-1})^i_k {\sf d} T^k_l (T^{-1})^l_j. 
\eqno(35)$$ 

Note that it is easy to check that $\hat{\epsilon}$ and $\hat{S}$ leave 
invariant the relations (15) and (16). Consequently, we can say that the 
structure $(\Gamma, \hat{\Delta}, \hat{\epsilon}, \hat{S})$ is a graded Hopf 
algebra. 

\noindent
{\bf C. Cartan-Maurer one-forms and their relations}

To complete the differential geometric scheme we need the Cartan-Maurer 
one-forms. As in analogy with the right-invariant one-forms on a Lie 
group in classical differential geometry, one can construct the matrix 
valued one-form $\Omega$ where 
$$\Omega = {\sf d} T ~T^{-1}. \eqno(36)$$
So we can write the matrix elements (right one-forms) of $\Omega$ as follows 
$$w_1 = \alpha A + b C, \qquad u = \alpha B + b D, $$
$$w_2 = \delta D + c B, \qquad v = c A + \delta C, \eqno(37)$$ 
where $T^{-1} = \left(\matrix{A&B\cr C&D}\right)$. We now wish to find 
the commutation relations of the matrix entries of $T$ with those of $\Omega$. 
So we need the commutation relations between the matrix elements of $T$ 
and $T^{-1}$, which may be computed directly, as follows: 
$$a A = A a + h (A - D) \beta, \qquad a B = B a, $$
$$a C = C a + h (1 - {\cal D}_h), \qquad a D = D a, $$
$$\beta A = A \beta, \qquad \beta B = - B \beta, $$
$$\beta C = - C \beta + h (D - A) \beta, \qquad \beta D = D \beta, $$
$$\gamma A = A \gamma + h (1 - {\cal D}_h^{-1}), \eqno(38)$$
$$\gamma B = - B \gamma + h (A - D) \beta, $$
$$\gamma C = - C \gamma, \qquad \gamma D = D \gamma + h ({\cal D}_h - 1), $$
$$d A = A d, \qquad 
  d C = C d + h ({\cal D}_h^{-1} - 1), $$
$$d B = B d, \qquad d D = D d + h (D - A)\beta. $$
Using these relations, we now find the commutation relations of the 
matrix entries of $T$ with those of $\Omega$ : 
$$a w_1 = w_1 a - h u a, \qquad a u = u a, $$
$$a v = v a + h(w_1 + w_2) a, \qquad a w_2 = w_2 a + h u a, $$
$$\beta w_1 = - w_1 \beta + h u \beta, \qquad \beta u = u \beta, $$
$$\beta v = v \beta + h (w_1 + w_2) \beta, \qquad 
  \beta w_2 = - w_2 \beta - h u \beta, \eqno(39)$$
$$\gamma w_1 = - w_1 \gamma + h (2 w_1 a + u \gamma), \qquad 
  \gamma u = u \gamma + 2 h u a, $$
$$\gamma v = v \gamma + h (w_1 \gamma + 2 v a + w_2 \gamma), \qquad 
  \gamma w_2 = - w_2 \gamma + h (2 w_2 a - u \gamma), $$
$$d w_1 = w_1 d - h (2 w_1 \beta + u d), \qquad 
  d u = u d + 2 h u \beta, $$
$$d v = v d + h (w_1 d + 2 v \beta + w_2 d), \qquad 
  d w_2 = w_2 d + h (u d - 2 w_2 \beta). $$

To obtain the commutation relations among the right Cartan-Maurer 
one-forms, we use the commutation relations of the matrix elements 
of $T^{-1}$ with the differentials of the group parameters 
which are given in the following: 
$$A \alpha = \alpha A + h (b A - \alpha B), \qquad A b = b A + h b B, $$
$$A c = c A - h (\alpha A + \delta A - c B), \qquad 
  A \delta = \delta A - h (b A + \delta B), $$
$$B \alpha = - \alpha B - h b B, \qquad B b = b B, $$
$$B c = c B - h (\alpha + \delta) B, \qquad 
  B \delta = - \delta B + h b B, $$
$$C \alpha = - \alpha C - h (\alpha A + \alpha D + b C), \qquad 
  C b = b C - h b (A + D), $$
$$C c = c C - h (\alpha C + c A + c D + \delta C), $$
$$C \delta = - \delta C + h (b C - \delta A - \delta D), \eqno(40)$$
$$D \alpha = \alpha D + h (\alpha B + b D), \qquad D b = b D - h b B, $$
$$D c = c D - h (\alpha D + c B + \delta D), \qquad 
  D \delta = \delta D + h (\delta B - b D).$$

Using these relations, we obtain the commutation relations of the right 
Cartan-Maurer forms with the differentials of the matrix elements of $T$ as 
follows:
$$w_1 \alpha = - \alpha w_1 - h \alpha u, \qquad 
  w_1 b = b w_1 - h b u, $$
$$w_1 c = c w_1 - h c u, \qquad 
  w_1 \delta = - \delta w_1 - h \delta u, $$
$$u \alpha = \alpha u, \qquad u b = b u, $$
$$u c = c u, \qquad u \delta = \delta u, \eqno(41)$$
$$v \alpha = \alpha v + h \alpha (w_1 - w_2), \qquad 
  v b = b v - h b (w_1 - w_2), $$
$$v c = c v - h c (w_1 - w_2), \qquad 
  v \delta = \delta v + h \delta (w_1 - w_2), $$
$$w_2 \alpha = - \alpha w_2 - h \alpha u, \qquad 
  w_2 b = b w_2 - h b u, $$
$$w_2 c = c w_2 - h c u, \qquad 
  w_2 \delta = - \delta w_2 - h \delta u.$$

We now obtain the commutation relations of the right Cartan-Maurer forms 
$$w_1 u = u w_1 - 2h u^2, \qquad w_2 u = u w_2, $$
$$w_1 v = v w_1 + 2h (w_1 w_2 - uv), \qquad w_2 v = v w_2, $$
$$w_1 w_2 = - w_2 w_1 - 2h u w_2, $$
$$w_1^2 = - 2h u w_1, \qquad w_2^2 = 0, \eqno(42)$$
$$u v = v u - 2h u w_2. $$

It can be checked that the elements ${\cal D}_h$ and $\hat{\cal D}$ 
commute with the Cartan-Maurer one-forms, i.e, both of the ${\cal D}_h$ 
and $\hat{\cal D}$ are still central elements. In the $q$-deformation, 
${\cal D}_q$ does not commute with the Cartan Maurer forms. 

Of course, the relations (4), (15), (16) and (38)-(42) can be obtained with 
the help of a matrix $R$ that acts on the square tensor space of the 
supergroup. The matrix $R$ is a solution of the quantum supergroup equation. 
The quantum supergroup relations (4) follows from the equation 
$$R T_1 T_2 = T_2 T_1 R,$$
where, in usual gradin tensor notation, 
$T_1 = T \otimes I$ and $T_2 = I \otimes T$ and 
$$R = \left(\matrix{1&0&0&0\cr -h&1&0&0\cr h&0&1&0\cr 0&h&h&1 \cr}\right).$$
The relations (15) are equivalent to the equation 
$$T'_1 \hat{T}_2 = R^{-1} \hat{T}_2 T_1 R,$$
where 
$$T'_1 = (-1)^{p(T_1)} T_1, \qquad \hat{T}_2 = {\sf d} T_2.$$
Applying the exterior differential {\sf d} on both sides of the above 
equation, one has 
$$(\hat{T}_1)' \hat{T}_2 = R \hat{T}_2' \hat{T}_1 R, \qquad 
  (\hat{T}_1)' = {\sf d} (T_1'),$$
which is equivalent to the relations (16). Similarly, the relations (39), (41) 
and (42 can be written, in a compact form, as follows, respectively: 
$$T'_1 \Omega_2 = R^{-1} \Omega_2 R T_1, $$
$$\hat{T}'_1 \Omega_2 = R \Omega_2' R \hat{T}_1, $$
$$\Omega_1' R^{-1} \Omega_2 R = - R \Omega_2' R \Omega_1. $$

Note that one can check that the action of {\sf d} on (39), (41) and also 
(42) is consistent. These relations allow us to evaluate the superalgebra of 
GL$_h(1\vert 1)$ by relating the generators of the superalgebra to the 
right one-forms. 

\noindent
{\bf IV. QUANTUM SUPERALGEBRA}

The commutation relations of Cartan-Maurer forms allow us to construct the 
algebra of the generators. To obtain the quantum superalgebra of the 
algebra generators we first write the Cartan-Maurer forms as 
$$\alpha = w_1 a + u \gamma, \qquad b = w_1 \beta + u d, $$
$$c = w_2 d + v \beta, \qquad \delta = w_2 \gamma + v a.   \eqno(43)$$
The differential {\sf d} can then the expressed in the form 
$${\sf d} = w_1 T_1 + w_2 T_2 + u \nabla_+ + v \nabla_-. \eqno(44)$$
Here $T_1$, $T_2$ and $\nabla_{\pm}$ are the quantum algebra generators. 
We now shall obtain the commutation relations of these generators. 
Considering an arbitrary function $f$ of the matrix elements of $T$ and using 
the nilpotency of the exterior differential {\sf d} one has 
$$({\sf d} w_i) T_i f + ({\sf d} u_i) \nabla_i f = 
  w_i {\sf d} T_i f - u_i {\sf d} \nabla_i f, \eqno(45)$$
where 
$$w_i \in \{w_1,w_2\}, \qquad u_i \in \{u,v\}, \qquad 
  \nabla_i \in \{\nabla_+, \nabla_-\}.$$
So we need the four two-forms. To obtain these, using the nilpotency of the 
differential {\sf d}, we can write ${\sf d} \Omega$ of the form 
$${\sf d} \Omega = \sigma_3 \Omega \sigma_3 \Omega, \qquad 
  \sigma_3 = \left(\matrix{1 & 0 \cr 0 & - 1 \cr}\right). \eqno(46)$$
In terms of the two-forms, these become 
$${\sf d} w_1 = w_1^2 - u v, \qquad {\sf d} u = w_1 u - u w_2, $$
$${\sf d} w_2 = w_2^2 - v u, \qquad {\sf d} v = w_2 v - v w_1. \eqno(47)$$
Using the Cartan-Maurer equations we find the following commutation 
relations for the quantum superalgebra: 
$$T_1 T_2 - T_2 T_1 = 2h \nabla_- T_1, $$
$$T_1 \nabla_+ - \nabla_+ T_1 = - \nabla_+ + 2h (T_1^2 - T_1), $$
$$T_2 \nabla_+ - \nabla_+ T_2 = \nabla_+ - 2h 
  (T_2 T_1 + T_2 - \nabla_+ \nabla_-), $$
$$T_1 \nabla_- - \nabla_- T_1 = \nabla_-, $$
$$T_2 \nabla_- - \nabla_- T_2 = - \nabla_-, \eqno(48)$$
$$\nabla_+^2 = - 2h T_1 \nabla_+, \qquad \nabla_-^2 = 0,$$
$$\nabla_- \nabla_+ + \nabla_+ \nabla_- = T_1 + T_2 - 2h \nabla_- T_1.$$

The commutation relations (48) of the algebra generators should be consistent 
with monomials of the matrix elements of $T$. To do this, we evaluate the 
commutation relations between the generators of algebra and the matrix 
elements of $T$. The commuation relations of the generators with the matrix 
elements can be extracted from the Leibniz rule: 
$${\sf d} (a f) = ({\sf d} a) f + a ({\sf d} f) \Longrightarrow 
  (w_i T_i + u_i \nabla_i) a = {\sf d} a + a (w_i T_i + u_i \nabla_i), 
  \eqno(49)$$
etc. This yields 
$$T_1 a  = a + a T_1 - h a \nabla_-, $$
$$T_1 \beta  = \beta + \beta T_1 + h \beta \nabla_-, $$
$$T_1 \gamma  = \gamma T_1 + h (2 a T_1 + \gamma \nabla_-), $$
$$T_1 d  = d T_1 + h (2 \beta T_1 - d \nabla_-), $$
$$T_2 a  = a T_2 - h a \nabla_-, $$
$$T_2 \beta  = \beta T_2 + h \beta \nabla_-, $$
$$T_2 \gamma  = \gamma + \gamma T_2 + h (2a T_2 + \gamma \nabla_-), $$
$$T_2 d = d + d T_2 + h (2 \beta T_2 - d \nabla_-), \eqno(50)$$
$$\nabla_+ a = \gamma + a \nabla_+ - h a (T_1 - T_2), $$
$$\nabla_+ \beta = d - \beta \nabla_+ - h \beta (T_1 - T_2), $$
$$\nabla_+ \gamma = - \gamma \nabla_+ - 
  h (2a \nabla_- + \gamma T_1 - \gamma T_2), $$
$$\nabla_+ d = d \nabla_+ + h (2\beta \nabla_+ - d T_1 + d T_2), $$
$$\nabla_- a = a \nabla_-,  \qquad \nabla_- \beta = - \beta \nabla_-, $$
$$\nabla_- \gamma = a -  \gamma \nabla_- - 2h a \nabla_-, $$
$$\nabla_- d = \beta + d \nabla_- + 2 h \beta \nabla_-. $$

\noindent
{\bf V. CONCLUSION}

To conclude, we introduce here commutation relations between the 
group parameters and their partial derivatives and thus illustrate 
the connection between the relations in sec. IV, and the relations 
which will be now obtained. 

To proceed, let us first obtain the relations of the group parameters 
with their partial derivatives. We know that the right exterior 
differential {\sf d} can be expressed in the form 
$${\sf d} f = (\alpha \partial_a + b \partial_\beta + c \partial_\gamma + 
  \delta \partial_d) f. \eqno(51)$$
Then, replacing $f$ with $af$, etc. we obtain the following commutation 
relations 
$$\partial_a a = 1 + a \partial_a - 
  h (\beta \partial_a + a \partial_\gamma), $$
$$\partial_a \beta = \beta \partial_a + h \beta \partial_\gamma, $$
$$\partial_a \gamma = \gamma \partial_a + h (a \partial_a + d \partial_a + 
  \gamma \partial_\gamma), $$
$$\partial_a d = d \partial_a + h (\beta \partial_a - d \partial_\gamma), $$
$$\partial_\beta a = a \partial_\beta - h (a \partial_a - a \partial_d 
  + \beta \partial_\beta), $$
$$\partial_\beta \beta = 1 - \beta \partial_\beta - 
   h \beta (\partial_a - \partial_d), $$
$$\partial_\beta \gamma = - \gamma \partial_\beta - 
  h (a \partial_\beta + \gamma \partial_a - \gamma \partial_d + 
   d \partial_\beta), $$
$$\partial_\beta d = d \partial_\beta + h (\beta \partial_\beta - d \partial_a 
  + d \partial_d), \eqno(52)$$
$$\partial_\gamma a = a \partial_\gamma -h \beta \partial_\gamma, 
  \qquad \partial_\gamma \beta = - \beta \partial_\gamma, $$
$$\partial_\gamma \gamma = 1 - \gamma \partial_\gamma - h 
  (a \partial_\gamma + d \partial_\gamma), $$
$$\partial_\gamma d = d \partial_\gamma + h \beta \partial_\gamma, $$
$$\partial_d a = a \partial_d - h (a \partial_\gamma + \beta \partial_d), $$
$$\partial_d \beta = \beta \partial_d + h \beta \partial_\gamma, $$
$$\partial_d \gamma = \gamma \partial_d + h (a \partial_d + 
  \gamma \partial_\gamma + d \partial_d), $$
$$\partial_d d = 1 + d \partial_d + h (\beta \partial_d - d \partial_\gamma).$$
We thus find the commutation relations between the derivatives. 
These relations can be obtained by using the nilpotency of the right exterior 
differential {\sf d} and they have the form 
$$\partial_a \partial_\beta = \partial_\beta \partial_a + h (
   \partial_d \partial_a - \partial_\beta \partial_\gamma - \partial_a^2), $$
$$\partial_d \partial_\beta = \partial_\beta \partial_d - h (
   \partial_a \partial_d + \partial_\beta \partial_\gamma - \partial_d^2), $$
$$\partial_a \partial_\gamma = \partial_\gamma \partial_a, \qquad 
  \partial_d \partial_\gamma = \partial_\gamma \partial_d, $$
$$\partial_\beta \partial_\gamma = - \partial_\gamma \partial_\beta + h 
   \partial_\gamma (\partial_a - \partial_d), $$
$$ \partial_\beta^2 = h \partial_\beta (\partial_a - \partial_d), \qquad 
   \partial_\gamma^2 = 0, \eqno(53)$$
$$ \partial_a \partial_d = \partial_d \partial_a + h \partial_\gamma (
  \partial_d - \partial_a). $$
The (graded) Hopf algebra structure for $\partial$ is given by 
$$\Delta(\partial_a) = \partial_a \otimes \partial_a + 
  \partial_\beta \otimes \partial_\gamma, \qquad 
  \Delta(\partial_\beta) = \partial_a \otimes \partial_\beta + 
  \partial_\beta \otimes \partial_d, $$
$$\Delta(\partial_d) = \partial_d \otimes \partial_d + 
  \partial_\gamma \otimes \partial_\beta, \qquad 
  \Delta(\partial_\gamma) = \partial_\gamma \otimes \partial_a + 
  \partial_d \otimes \partial_\gamma, \eqno(54)$$
$$\varepsilon(\partial_a) = 1 = \varepsilon(\partial_d), \qquad 
  \varepsilon(\partial_\beta) = 0 = \varepsilon(\partial_\gamma), $$
$$S(\partial_a) = \partial_a^{-1} + \partial_a^{-1}  \partial_\beta 
   \partial_d^{-1} \partial_\gamma \partial_a^{-1}, \qquad 
  S(\partial_\beta) = - \partial_a^{-1} \partial_\beta \partial_d^{-1}, $$
$$S(\partial_d) = \partial_d^{-1} + \partial_d^{-1}  \partial_\gamma 
   \partial_a^{-1} \partial_\beta \partial_d^{-1}, \qquad 
  S(\partial_\gamma) = - \partial_d^{-1} \partial_\gamma \partial_a^{-1}, $$
provided that the formal inverses $\partial_a^{-1}$ and $\partial_d^{-1}$ 
exist. However these co-maps do not leave invariant the relations (52). 

We know, from Sec. IV, that the right exterior differential {\sf d} can 
be expressed in the form (31), which we repeat here, 
$${\sf d} f = (w_1 T_1 + u \nabla_+ + v \nabla_- + w_2 T_2) f. \eqno(55)$$
Considering (51) together (55) and using (43) one has 
$$T_1 = a \partial_a + \beta \partial_\beta, \qquad 
  \nabla_+ = \gamma \partial_a + d \partial_\beta, $$
$$T_2 = d \partial_d + \gamma \partial_\gamma, \qquad 
  \nabla_- = a \partial_\gamma + \beta \partial_d. \eqno(56)$$
Using the relations (52) and (53) one can check that the relations of the 
generators in (56) coincide with (48). It can also be verified that, 
the action of the generators in (56) on the group parameters coincide 
with (50). 

The classical limit $h \longrightarrow 0$ of the right-invariant differential 
calculus is the undeformed (ordinary) differential calculus. 

Note that if we make the identification 
$$u \longrightarrow {x\over 2}, \qquad w_1 \longrightarrow \theta,$$
where $x$ and $\theta$ are the coordinates of superplane, we have 
$$x \theta = \theta x + h x^2, \theta^2 = - h x \theta.$$
One of the interesting problems may be to construct linear connections$^{15}$ 
on the $h$-superplane. 

\noindent
{\bf ACKNOWLEDGEMENT}

This work was supported in part by {\bf T. B. T. A. K.} the 
Turkish Scientific and Technical Research Council. 

\noindent
$^1$ N. Reshetikhin, L. Takhtajan and L. Faddeev, 
     Leningrad Math. J. {\bf 1}, 193 (1990); 
     S. Majid, Int. J. Mod. Phys. A {\bf 5}, 1 (1990).\\
$^2$ Yu I. Manin, Commun. Math. Phys. {\bf 123}, 163 (1989); 
     E. Corrigan, D. Fairlie, P. Fletcher and R. Sasaki, 
     J. Math. Phys. {\bf 31}, 776 (1990); 
     W. Schmidke, S. Vokos and B. Zumino, Z. Phys. C {\bf 48}, 249 (1990).\\ 
$^3$ L. Alvarez-Gaume, C. Gomes and G.Sierra, 
     Nucl. Phys. B {\bf 319}, 155 (1989); 
     T. Curtright, D. Fairlie and C. Zachos, "Quantum groups", in Proc. 
     Argonne Workshop (World Scientific, 1990); 
     D. Fairlie and C. Zachos, Phys. Lett. B {\bf 256}, 43 (1991).\\
$^4$ J. Madore, {\it An introduction to Noncommutative Geometry and its 
     Physical Applications} (Cambridge U. P., Cambridge, 1995).
$^5$ S. L. Woronowicz, Commun. Math. Phys. {\bf 122}, 125 (1989).\\
$^6$ S. L: Woronowicz, Kyoto Univ. {\bf 23}, 117 (1987). \\
$^7$ A. Schirmacher, J. Wess and B. Zumino, Z. Phys. C {\bf 49}, 317 (1990). 
$^8$ P. Aschieri and L. Castellani, Int. J. Mod. Phys. A {\bf 8}, 1667 (1993); 
     B. Jurco, Lett. Math. Phys. {\bf 22}, 177 (1991); 
     A. Sudbery, Phys. Lett. B {\bf 284}, 61 (1992); 
     F. Muller-Hoissen, J. Phys. A {\bf 25}, 1703 (1992).\\
$^9$ S. Zakrzewski, Lett. Math. Phys. {\bf 22}, 287 (1991); 
     B. A. Kupershmidt, J. Phys. A {\bf 25}, L1239 (1992); 
     Ch. Ohn, Lett. Math. Phys. {\bf 25}, 89 (1992). \\
$^{10}$ A. Aghamohammadi, M. Khorrami and A. Shariati, 
        J. Phys. A {\bf 28}, L225 (1995).\\
$^{11}$ V. Karimipour, Lett. Math. Phys. {\bf 25}, 87 (1994); 
        {\bf 35}, 303 (1995). \\
$^{12}$ L. Dabrowski and P. Parashar, Lett. Math. Phys. {\bf 38}, 331 (1996).\\
$^{13}$ S. Celik and S. A. Celik, J. Phys. A {\bf 31}, 9685 (1998). \\
$^{14}$ S. Celik, Balkan Phys. Lett. {\bf 5}, 149 (1997).\\ 
$^{15}$ Y. Georgelin, T. Masson, and J.-C. Wallet, "Linear Connections on the 
        Two-parameter Quantum Plane", q-alg/9507032 (1995). 

\end{document}